\documentclass{amsart} 
\usepackage{latexsym} 
\usepackage{amssymb} 
\usepackage{amsmath}
\usepackage{enumerate}
\usepackage{hyperref}  

\begin{document}


\newtheorem{question}{Question}  
\newtheorem{theorem}{Theorem} 
\newtheorem{problem}{Problem} 
\newtheorem{definition}{Definition} 
\newtheorem{lemma}{Lemma} 
\newtheorem{proposition}{Proposition} 
\newtheorem{corollary}{Corollary} 
\newtheorem{example}{Example} 
\newtheorem{conjecture}{Conjecture} 
\newtheorem{algorithm}{Algorithm} 
\newtheorem{exercise}{Exercise} 
\newtheorem{remarkk}{Remark} 
 
\newcommand{\be}{\begin{equation}} 
\newcommand{\ee}{\end{equation}} 
\newcommand{\bea}{\begin{eqnarray}} 
\newcommand{\eea}{\end{eqnarray}} 

\newcommand{\eeq}{\end{equation}} 

\newcommand{\eeqn}{\end{eqnarray}} 
\newcommand{\beaa}{\begin{eqnarray*}} 
\newcommand{\eeaa}{\end{eqnarray*}} 

\newcommand{\lip}{\langle} 
\newcommand{\rip}{\rangle}

\newcommand{\uu}{\underline} 
\newcommand{\oo}{\overline} 
\newcommand{\La}{\Lambda} 
\newcommand{\la}{\lambda} 
\newcommand{\eps}{\varepsilon} 
\newcommand{\om}{\omega} 
\newcommand{\Om}{\Omega} 
\newcommand{\ga}{\gamma} 
\newcommand{\rrr}{{\Bigr )}} 
\newcommand{\qqq}{{\Bigl\|}} 
 
\newcommand{\dint}{\displaystyle\int} 
\newcommand{\dsum}{\displaystyle\sum} 
\newcommand{\dfr}{\displaystyle\frac} 
\newcommand{\bige}{\mbox{\Large\it e}} 
\newcommand{\integers}{{\Bbb Z}} 
\newcommand{\rationals}{{\Bbb Q}} 
\newcommand{\reals}{{\rm I\!R}} 
\newcommand{\realsd}{\reals^d} 
\newcommand{\realsn}{\reals^n} 
\newcommand{\NN}{{\rm I\!N}} 
\newcommand{\DD}{{\rm I\!D}} 
\newcommand{\degree}{{\scriptscriptstyle \circ }} 
\newcommand{\dfn}{\stackrel{\triangle}{=}} 
\def\complex{\mathop{\raise .45ex\hbox{${\bf\scriptstyle{|}}$} 
     \kern -0.40em {\rm \textstyle{C}}}\nolimits} 
\def\hilbert{\mathop{\raise .21ex\hbox{$\bigcirc$}}\kern -1.005em {\rm\textstyle{H}}} 
\newcommand{\RAISE}{{\:\raisebox{.6ex}{$\scriptstyle{>}$}\raisebox{-.3ex} 
           {$\scriptstyle{\!\!\!\!\!<}\:$}}} 
 
\newcommand{\hh}{{\:\raisebox{1.8ex}{$\scriptstyle{\degree}$}\raisebox{.0ex} 
           {$\textstyle{\!\!\!\! H}$}}} 

\newcommand{\OO}{\won} 
\newcommand{\calA}{{\mathcal A}} 
\newcommand{\calB}{{\cal B}} 
\newcommand{\calC}{{\cal C}} 
\newcommand{\calD}{{\cal D}} 
\newcommand{\calE}{{\cal E}} 
\newcommand{\calF}{{\mathcal F}} 
\newcommand{\calG}{{\cal G}} 
\newcommand{\calH}{{\cal H}} 
\newcommand{\calK}{{\cal K}} 
\newcommand{\calL}{{\mathcal L}} 
\newcommand{\calM}{{\mathcal M}} 
\newcommand{\calO}{{\cal O}} 
\newcommand{\calP}{{\cal P}} 
\newcommand{\calU}{{\mathcal U}} 
\newcommand{\calX}{{\cal X}} 
\newcommand{\calXX}{{\cal X\mbox{\raisebox{.3ex}{$\!\!\!\!\!-$}}}} 
\newcommand{\calXXX}{{\cal X\!\!\!\!\!-}} 
\newcommand{\gi}{{\raisebox{.0ex}{$\scriptscriptstyle{\cal X}$} 
\raisebox{.1ex} {$\scriptstyle{\!\!\!\!-}\:$}}} 
\newcommand{\intsim}{\int_0^1\!\!\!\!\!\!\!\!\!\sim} 
\newcommand{\intsimt}{\int_0^t\!\!\!\!\!\!\!\!\!\sim} 
\newcommand{\pp}{{\partial}} 
\newcommand{\al}{{\alpha}} 
\newcommand{\sB}{{\cal B}} 
\newcommand{\sL}{{\cal L}} 
\newcommand{\sF}{{\cal F}} 
\newcommand{\sE}{{\cal E}} 
\newcommand{\sX}{{\cal X}} 
\newcommand{\R}{{\rm I\!R}} 
\renewcommand{\L}{{\rm I\!L}} 
\newcommand{\vp}{\varphi} 
\newcommand{\N}{{\rm I\!N}} 
\def\ooo{\lip} 
\def\ccc{\rip} 
\newcommand{\ot}{\hat\otimes} 
\newcommand{\rP}{{\Bbb P}} 
\newcommand{\bfcdot}{{\mbox{\boldmath$\cdot$}}} 
 
\renewcommand{\varrho}{{\ell}} 
\newcommand{\dett}{{\textstyle{\det_2}}} 
\newcommand{\sign}{{\mbox{\rm sign}}} 
\newcommand{\TE}{{\rm TE}} 
\newcommand{\TA}{{\rm TA}} 
\newcommand{\E}{{\rm E\, }} 
\newcommand{\won}{{\mbox{\bf 1}}} 
\newcommand{\Lebn}{{\rm Leb}_n} 
\newcommand{\Prob}{{\rm Prob\, }} 
\newcommand{\sinc}{{\rm sinc\, }} 
\newcommand{\ctg}{{\rm ctg\, }} 
\newcommand{\loc}{{\rm loc}} 
\newcommand{\trace}{{\, \, \rm trace\, \, }} 
\newcommand{\Dom}{{\rm Dom}} 
\newcommand{\ifff}{\mbox{\ if and only if\ }} 
\newcommand{\nproof}{\noindent {\bf Proof:\ }} 
\newcommand{\nproofYWN}{\noindent {\bf Proof of Theorem~\ref{YWN}:\ }} 
\newcommand{\remark}{\noindent {\bf Remark:\ }} 
\newcommand{\remarks}{\noindent {\bf Remarks:\ }} 
\newcommand{\note}{\noindent {\bf Note:\ }}

\newcommand{\boldx}{{\bf x}} 
\newcommand{\boldX}{{\bf X}} 
\newcommand{\boldy}{{\bf y}} 
\newcommand{\boldR}{{\bf R}} 
\newcommand{\uux}{\uu{x}} 
\newcommand{\uuY}{\uu{Y}} 
 
\newcommand{\limn}{\lim_{n \rightarrow \infty}} 
\newcommand{\limN}{\lim_{N \rightarrow \infty}} 
\newcommand{\limr}{\lim_{r \rightarrow \infty}} 
\newcommand{\limd}{\lim_{\delta \rightarrow \infty}} 
\newcommand{\limM}{\lim_{M \rightarrow \infty}} 
\newcommand{\limsupn}{\limsup_{n \rightarrow \infty}} 
 
\newcommand{\ra}{ \rightarrow } 

 \newcommand{\mlim}{\lim_{m \rightarrow \infty}}  
 \newcommand{\limm}{\lim_{m \rightarrow \infty}}  
 \newcommand{\nlim}{\lim_{n \rightarrow \infty}} 
 
 
 
 
 
 
 
\newcommand{\one}{\frac{1}{n}\:} 
\newcommand{\half}{\frac{1}{2}\:} 
 
\def\le{\leq} 
\def\ge{\geq} 
\def\lt{<} 
\def\gt{>} 
 
\def\squarebox#1{\hbox to #1{\hfill\vbox to #1{\vfill}}} 
\newcommand{\nqed}{\hspace*{\fill} 
           \vbox{\hrule\hbox{\vrule\squarebox{.667em}\vrule}\hrule}\bigskip} 

 
  \author{R\'{e}mi  Lassalle}
\address[R\'{e}mi  Lassalle]{CEREMADE Universit\'{e} Paris-Dauphine Place du Mar\'{e}chal de Lattre de Tassigny 75016 Paris }
\email{lassalle@ceremade.dauphine.fr}

\author{Ana Bela Cruzeiro}
\address[Ana Bela Cruzeiro]{GFMUL and Departamento Matem\'{a}tica Instituto Superior T\'{e}cnico, Univ. de Lisboa, Av. Rovisco Pais, 1049-001 Lisbon, Portugal }
\email{abcruz@math.tecnico.ulisboa.pt }

\title[Symmetries and martingales]{Symmetries and martingales in a stochastic model for  the Navier-Stokes equation}

\maketitle 
\noindent 
{\bf Abstract: }{\small{A stochastic  description of solutions of the Navier-Stokes equation is investigated. These solutions are represented by laws of finite dimensional semi-martingales and characterized by a weak Euler-Lagrange condition. A least action principle,  related to the relative entropy, is provided. Within this stochastic framework, 
by assuming further symmetries, the corresponding invariances are expressed by martingales, stemming from a weak Noether's theorem.}
\\ 

\vspace{0.5cm}

 
\noindent 
\tableofcontents


Several stochastic models for the \textit{Navier-Stokes} equation have been proposed in the literature. Some refer to random perturbations of the fluid velocity. This is not the case here: we are interested in stochastic Lagrangian paths whose (mean) velocity, or drift, represent the deterministic velocity of the fluid. Different studies  of stochastic Lagrangian paths in fluid dynamics and in particular in turbulence can be found in a collection of works, from which we refer to \cite{BJT}, \cite{MR1} and \cite{MR2} as examples.
Also representation formulae in terms of different random processes
were given in \cite{C}, \cite{CI}, \cite{CRUZSH}, among others. 

Concerning the derivation of solutions of Navier-Stokes equations from (stochastic)
variational principles, after the early articles \cite{NYZ} and \cite{Y}, such principles were developed in \cite{CRUZCIP} and subsequent works. We
mention also \cite{IF} and \cite{G} for different,  unrelated  approaches to the same kind of problems.

 In  \cite{CRUZLAS} and \cite{RLZ} a weak description of a stochastic deformation of mechanics has been investigated: the Euler-Lagrange condition extends as a condition on laws of stochastic processes. The main originality of this approach, with respect to other weak deformations of mechanics, is to handle problems in a functional-analytic framework, in  a full consistency with the deterministic case.

In this paper, within the framework of \cite{CRUZLAS} and \cite{RLZ}, we develop a specific case, from an example of \cite{CRUZLAS}, which provides a  \textit{stochastic description} for solutions to the \textit{Navier-Stokes equation}.  Then we investigate several \textit{symmetries}, whose associated invariances correspond in this setup to  \textit{martingales}.  

 Section~\ref{1} fixes the framework and notations of the  paper; the \textit{weak Euler-Lagrange condition} of \cite{RLZ} and \cite{CRUZLAS} is recalled. Under 
 conditions, in Section~\ref{2}, a map $$P : u\to P_u ,$$  associates laws of $\mathbb{R}^d-$valued \textit{semi-martingales} to \textit{divergence free vector fields}. 
 Solutions of the \textit{Navier-Stokes equation}   are shown to be  divergence free vector fields $u$, whose associated probability $P_u$ satisfies a weak \textit{Euler-Lagrange condition} (Proposition~\ref{ELP}); Corollary~\ref{LAPNS} characterizes those solutions as \textit{critical points} of the \textit{stochastic action} $$\mathcal{S}^p (\nu) := E_\nu\left[\int_0^1 \left(\frac{|v_s^\nu|^2}{2} -p(1-t,W_t) \right) dt\right],$$ where $\nu$ denotes the \textit{law} of specific continuous \textit{semi-martingales}, and  $(v_t^\nu)$ the characteristic drift of $\nu$, as stated accurately below. The function $p$ is a smooth  \textit{pressure field} which is assumed to be \textit{given}. 
 
 The action functional above  is related to the relative \textit{entropy} with respect to a reference law $\mu_p$ induced by the pressure field.  
 
 Finally, within this  stochastic model, Section~\ref{3} investigates \textit{invariances}, stemming from \textit{symmetries}, by the weak \textit{Noether's theorem} of \cite{CRUZLAS}; within this stochastic framework  \textit{martingales} on the canonical space play the r\^ole of constants of motion in classical mechanics.

\section{The weak stochastic Euler-Lagrange condition.}
\label{1}
The \textit{weak} stochastic \textit{Euler-Lagrange condition}, recalled below, was introduced in \cite{CRUZLAS},\cite{RLZ}. It \textit{embeds} in probability measures, specifically in a set of \textit{laws} of semi-martingales, the classical condition. Thus, it provides a \textit{functional analytic} approach to tackle \textit{stochastic variational problems}; in particular, in contrast with usual diffusion approaches, it directly embeds the  not stochastic case. In this context the extension of \textit{Noether's theorem} becomes natural.  Moreover,  as stated in \cite{CRUZLAS}, another specificity of this framework is that it provides critical conditions to \textit{semi-martingale optimal transportation problems}. The latter, introduced in \cite{Touzi}, correspond to a relaxation of a specific  \textit{dynamical Schr\"{o}dinger problems} (see \cite{LEO}) by allowing, in particular,  the characteristic dispersion to be not-trivial. Finally, as it is expected of  optimization over a subset of Borel probabilities on a Polish space, one crucial advantage  of this framework is that \textit{compactness} is rather simple to obtain.  

\subsection{Admissible trajectories.} 

Trajectories of infinitely small passive tracers in fluids can be described by elements in the space $W:=$ $C([0,1],\mathbb{R}^d)$, namely the set of continuous $\mathbb{R}^d$-valued paths  (where we consider the norm $|.|_W$ of uniform convergence), endowed with the Borel sigma-field $\mathcal{B}(W)$. In particular, trajectories of finite energy can be described by the subset $$H:=\left\{ h\in W , h:= \int_0^. \dot{h}_s ds, \int_0^1 |\dot{h}_s|_{\mathbb{R}^d}^2 ds <\infty\right\},$$ of absolutely continuous paths  with square integrable derivatives.  

Let $(W_t)_{t\in[0,1]}$  denote the \textit{evaluation process} $$(t,\omega) \in[0,1]\times W\to W_t(\omega):= \omega(t) \in \mathbb{R}^d.$$
We consider  $(\mathcal{F}_t^0)$ the natural (past) filtration and  denote by $\mathcal{P}_W$  the set of Borel probabilities on $W$.

\subsection{A weak description of random trajectories.} In order    to avoid measurability issues, we consider $(\mathcal{F}_t^\nu)$, the $\nu-$usual augmentation of the filtration   $(\mathcal{F}_t^0)$ under $\nu\in \mathcal{P}_W$. The latter naturally models random trajectories, since any $\nu\in \mathcal{P}_W$ is the law of the evaluation process, on the completed probability space $(W,\mathcal{B}(W)^\nu,\nu)$. We define $\mathbb{S}$ to be the subset  of $\nu\in \mathcal{P}_W$ such that there exists a $(\mathcal{F}
_t^\nu)-$martingale $(M^\nu_t)$ on $(W,\mathcal{B}(W)^\nu,\nu)$, which satisfies $$W_t=W_0 + M_t^\nu  + \int_0^t v_s^\nu ds, $$  for all $t\in [0,1]$, $\nu-a.s.$, where $(v_s^\nu)$ is a $(\mathcal{F}_t^\nu)-$\textit{predictable process} on the same space, and where the \textit{predicable covariation process} of $(M^\nu)$ is of the specific form $$<(M_t^\nu)^i, (M^\nu_t)^j> = \int_0^.(\alpha_s^\nu)^{ij} ds,$$  for  a predictable process $(\alpha_s^\nu)$; subsequently, by abuse of language, we refer to $(v_t^\nu,\alpha_t^\nu)$ as the \textit{characteristics} of $\nu$. In the whole paper, notations are those of \cite{CRUZLAS}.

\subsection{A weak Euler-Lagrange condition} Given a \textit{smooth Lagrangian} function  \begin{equation} \label{lagrangiandef} 
\mathcal{L} :(t,x, v, a)\in [0,1] \times \mathbb{R}^d\times \mathbb{R}^d \times (\mathbb{R}^d\otimes \mathbb{R}^d) \to \mathcal{L}_t(x,v,a)\in \mathbb{R}, \end{equation} the classical  \textit{Euler-Lagrange condition} naturally extends to $\mathbb{S}$ (see \cite{RLZ}).   A semi-martingale  $\nu\in \mathbb{S}$ satisfies the Euler-Lagrange condition  if there exists a $(\mathcal{F}_t^\nu)$ \textit{c\`{a}d-l\`{a}g martingale} $(N_t^\nu)$, such that \begin{equation} \label{eulerlagrange} \partial_{v} \mathcal{L}_t(W_t,v_t^\nu,\alpha_t^\nu) - \int_0^t \partial_q \mathcal{L}_s(W_s,v_s^\nu,\alpha_s^\nu) ds = N_t^\nu \quad  \lambda\otimes \nu-a.e., \end{equation}   $\partial_q \mathcal{L}$ and  $\partial_{v} \mathcal{L}_t$ denoting the respective gradients (in the first and in the second variables, respectively).

 \begin{remarkk} Similar conditions were considered in \cite{BISMUT} for arbitrary semi-martingales $U$ on abstract stochastic basis $(\Omega,\mathcal{A}, (\mathcal{A}_t), \mathcal{P})$.  On the contrary, condition~(\ref{eulerlagrange}) imposes  constraints on laws of processes. A semi-martingale $U$ on an arbitrary stochastic basis whose law satisfies~(\ref{eulerlagrange})  exhibits very precise properties, depending on the Lagrangian; for instance, in a specific case, it is associated to systems of coupled stochastic differential equations; the latter are not satisfied, in general,  when $U$ verifies the critical condition of \cite{BISMUT}. Moreover, the associated variational principles of \cite{BISMUT} do not contain the optimum criticality for semi-martingale optimal transportation problems either.  \end{remarkk}

\section{Navier-Stokes equation and the weak  Euler-Lagrange condition}  \label{2}

Henceforth, and \textit{until the end of the paper,  $p$ denotes a smooth map \begin{equation} \label{pression} p : (t,x)\in [0,1] \times \mathbb{R}^d \to  p(t,x) \in \mathbb{R}^+, \end{equation} which is further assumed to be bounded, with bounded derivatives}; $p$ models the \textit{pressure field}.  This map being given, we provide a \textit{stochastic model} for solutions of the equation \begin{equation} \label{NSp} \partial_t u + (u.\nabla) u =-\nabla p + \frac{\Delta u}{2} \ ;  \ div ~ u =0 .\end{equation} For  sake of clarity, we focus on the case where the \textit{divergence free velocity vector field}, involved in the Navier-Stokes equation, belongs to $$C_{b, div}^{1,2}([0,1]\times \mathbb{R}^d):= \{ u \in C^{1,2}([0,1]\times \mathbb{R}^d ; \mathbb{R}^d) \cap C_b([0,1]\times \mathbb{R}^d ; \mathbb{R}^d) :  div~ u(t,.)=0, \ for \ all \ t \in [0,1] \}.$$ 

\subsection{Description of dissipative flows by laws of semi-martingales}  Given $u \in C_{b, div}^{1,2}([0,1]\times \mathbb{R}^d)$,  we define $P_u$ to be the  probability measure, which is equivalent with respect to the \textit{Wiener measure $\mu\in \mathcal{P}_W$} (the law of \textit{standard Brownian motion}), with density  defined by $$\frac{dP_u}{d\mu}:= \exp\left(-\int_0^1 u(1-t, W_t) dW_t -\frac{1}{2}\int_0^1 |u(1-t,W_t)|^2 dt \right).$$ By the \textit{Girsanov theorem} (see for example  \cite{I-W}), we obtain a map \begin{equation} \label{Pdef} P : u \in C_{b, div}^{1,2}([0,1]\times \mathbb{R}^d) \to P_u \in \mathbb{S}, \end{equation} such that $\lambda\otimes P_u \ a.e.,$ $\alpha_s^{P_u}= I_{\mathbb{R}^d}$, and $v_t^{P_u}= -u(1-t,W_t),$ $\lambda$ denoting the Lebesgue measure.   \textit{It\^{o}'s formula} yields the following:
\begin{proposition}
\label{ELP}
A  time-dependent vector field $u\in C_{b,  div}^{1,2}([0,1]\times \mathbb{R}^d) $ satisfies the Navier-Stokes equation~(\ref{NSp}) if and only if $P_u$ satisfies~(\ref{eulerlagrange}) for \begin{equation} \label{lpdef} \mathcal{L}^p_t(x,v,a):= \frac{|v|^2}{2}- p(1-t,x).\end{equation}  \nqed \end{proposition}

\subsection{Stochastic action and relative entropy} Define the \textit{stochastic action} 
  associated to the Lagrangian $\mathcal{L}^p$ of~(\ref{lpdef}) by \begin{equation}\label{sp} \mathcal{S}^p : \nu \in \mathbb{S} \to  \mathcal{S}^p(\nu):= E_\nu\left[ \int_0^1\mathcal{L}^p_s(W_s, v_s^\nu, \alpha_s^\nu)ds \right] \in \mathbb{R} \cup\{+\infty\}. \end{equation}  As $\mathcal{L}^p$ and $P_u$ satisfy  the assumptions of Theorem 5.1 of \cite{CRUZLAS}, for all $p$ as above and $u\in C_{b,div}^{1,2}([0,1]\times \mathbb{R}^d)$, the $\mathbb{S}-$functional $\mathcal{S}^p$  is differentiable in the sense considered in \cite{CRUZLAS}. Thus, we obtain the following result :
\begin{proposition}
\label{LAPNS}
A vector field  $u\in C_{b,div}^{1,2}([0,1]\times \mathbb{R}^d)$ is a solution 
of the equation ~(\ref{NSp}) if and only if  $$\delta \mathcal{S}_{P_u}[h] =0,$$ for all 
$(\mathcal{F}_t^{P_u})-$adapted process $(h_t)$ of finite energy, such that  $$h_0=h_1=0 \ P_u-a.s.,$$ where  $P_u \in \mathbb{S}$ is given by~(\ref{Pdef}) and
$\delta \mathcal{S}_{P_u}$ denotes the $\mathbb{S}-$differential of \cite{CRUZLAS} at $P_u$. \nqed
\end{proposition}


\subsection{Least action principle and relative entropy.}
Subsequently, assuming $\nu\in \mathcal{P}_W$ to be  absolutely continuous with respect to a reference law $\eta\in \mathcal{P}_W$, 
 the relative \textit{entropy} of  $\nu$ w.r.t.  $\eta$ is defined as 
 $$\mathcal{H}(\nu | \eta) := E_\nu\left[\ln \frac{d\nu}{d\eta}\right].$$  By the representation formula in \cite{F1}, we obtain \begin{equation} 
\mathcal{S}^p(P_u) = \mathcal{H}(P_u | \mu_p) + \ln Z_p, \end{equation} for all 
$u\in C^{1,2}([0,1]\times \mathbb{R}^d),$ where $P_u$ is defined by~
(\ref{Pdef}), and where $\mu_p$ is the absolutely continuous probability with respect to the \textit{standard Wiener measure} ($W_0=0$, $\mu-a.s.$), whose density is given 
by $$\frac{d\mu_p}{d\mu}:= \frac{\exp\left(\int_0^1 p(1-s,W_s) ds
\right)}{Z_p},$$ 
with $Z_p$  a normalization constant. Whence we obtain the following ersatz of Proposition~\ref{LAPNS}: 
\begin{proposition}
\label{LAPNS2} A time-dependent vector field 
 $u\in C_{b,div}^{1,2}([0,1]\times \mathbb{R}^d)$ is  a solution 
of the equation ~(\ref{NSp})  if and only if $$\delta \mathcal{H}(. | \mu_p)_{P_u}[h] =0$$ for all 
$(\mathcal{F}_t^\nu)-$adapted process $(h_t)$ of finite energy, such that  $$h_0=h_1=0 ~~\ P_u-a.s..$$ \nqed
\end{proposition}

\section{Invariances and the stochastic Noether theorem.} \label{3}
Within this section we consider the case $d=3$ and we denote by $(e_{1},e_{2},e_{3})$ the canonical orthogonal basis of $\mathbb{R}^3$.
 We further assume that  $u \in C_{b, div}^{1,2}([0,1]\times \mathbb{R}^3)$ is a solution of~(\ref{NSp}) for a given smooth function $p$.
  By Proposition~\ref{ELP}, the associated {law} of the continuous {semi-martingale} $P_u$ (c.f. ~(\ref{Pdef})) satisfies~(\ref{eulerlagrange}) for $\mathcal{L}^p$ defined by~(\ref{lpdef}). 
  
  The next subsections investigate different \textit{symmetries} and compute the related local \textit{martingales}, stemming from the \textit{weak Noether Theorem} 6.1. of \cite{CRUZLAS}. In each particular case considered below \textit{the symmetries are expressed through a condition on the pressure field $p$}. Given  the associated family of transformations $(h^\epsilon)$, subsequently, the \textit{symmetry condition} on $p$ yields that $(h^\epsilon)$ is a smooth family of $\mathbb{S}-$invariant transformations for $\mathcal{L}^p$, in the sense considered in \cite{CRUZLAS}. 
  
   We recall this  \textit{symmetry condition} on $\mathbb{S}$. First, by setting  $$\Gamma^\epsilon : \omega \in W \to \Gamma^\epsilon(\omega) \in W,$$ where $$\Gamma^\epsilon_t(\omega):= h^\epsilon(t,\omega(t)),$$ for all $t\in [0,1]$, $\omega \in W$,  $(h^\epsilon)$ induces a family $(\Gamma^\epsilon)$ of transformations of $W$. Given $\eta\in \mathbb{S}$, for all $\epsilon$, $(\Gamma^\epsilon_t)$ defines a stochastic process on the probability space $(W,\mathcal{B}(W)^\eta, \eta)$. Thus, by \textit{It\^{o}'s formula} on the probability space $(W,\mathcal{B}(W)^\eta, \eta)$,  for all $\epsilon \in \mathbb{R}$, \textit{the transformation $h^\epsilon$ of the state space $\mathbb{R}^3$ is lifted to a transformation  $$\eta \in \mathbb{S} \to \Gamma^\epsilon_\star \eta \in \mathbb{S}$$ of $\mathbb{S}$}, by pushforward. The \textit{symmetry condition} considered in  \cite{CRUZLAS} consists in the relation  \begin{equation} \label{invcond} \mathcal{L}_t^p(W_t, v_t^\eta, \alpha_t^\eta) = \mathcal{L}_t^p(\Gamma^\epsilon_t, v_t^{\Gamma^\epsilon_\star \eta} \circ\Gamma^{\epsilon}, \alpha_t^{\Gamma^\epsilon_\star \eta}\circ \Gamma^\epsilon), \end{equation} 
holding a.e., for all $\eta \in \mathbb{S}$ in the domain of the map defined in~(\ref{sp}), and for all $\epsilon \in \mathbb{R}$. Here $\circ$ denotes the pullback of the ($\Gamma^\epsilon_\star \eta-$equivalence class of) map(s) $v_t^{\Gamma^\epsilon_\star \eta} : W\to \mathbb{R}^3$ with the ($\eta-$ equivalence class of) map(s) $\Gamma^\epsilon : W \to W$.

Consider  a smooth Lagrangian  $\mathcal{L}$ and assume that $\nu \in \mathbb{S}$ satisfies the weak Euler-Lagrange condition for this Lagrangian. The stochastic weak Noether's Theorem in  \cite{CRUZLAS}  associates to a family $(h_\epsilon )$ of $\mathbb{S} $-invariant transformations  of $ \mathcal{L}$ local martingales on the probability space $(W, \mathcal{B}(W)^\nu,\nu)$. These local martingales, that we denote by $\left(\mathcal{I}_t \right)_{t\in[0,1)}$ , are in fact explicitly given:

\begin{equation}\mathcal{I}_t:= <\frac{d}{d\epsilon}|_{\epsilon=0}h^\epsilon_t(W_t), \phi_t^\nu >_{\mathbb{R}^d} - \sum_i \left[ \frac{d}{d\epsilon}|_{\epsilon=0}{h^\epsilon_.(W_.)}^i, {\phi_.^\nu}^i \right]_t + \int_0^t  \theta_s ds \end{equation} where $[. ,. ]$ stands for the quadratic co-variation process of \textit{c\`{a}d-l\`{a}g} semi-martingales, $(\phi_t^\nu)$ denotes a \textit{c\`{a}dl\`{a}g} modification of the process $\partial_v \mathcal{L}_t(W_t,v_t^\nu, \alpha_t^\nu)$,
and
\begin{equation}\theta_s :=  \sum_{i,j}{{\kappa}_s^{i,j}} \frac{\partial \mathcal{L}}{\partial \alpha_{i,j}}(W_s,v_s^\nu, \alpha_s^\nu), \end{equation} 
 where  $(\kappa_s(\omega))$ is the $\mathcal{M}_d(\mathbb{R})-$valued process defined by  \begin{equation} \label{kappadl} {\kappa}_s(\omega):=  {\alpha^\nu_s}.
\left(\left(\nabla \frac{d}{d\epsilon}  h^\epsilon|_{\epsilon=0}\right)(s,W_s) \right)^\dagger +  \left(\left(\nabla \frac{d}{d\epsilon}  h^\epsilon|_{\epsilon=0}\right)(s,W_s)\right)  {\alpha^
\nu_s}. \end{equation}

\subsection{Symmetry by translation and the momentum process} Assume the \textit{symmetry by translation} along $e_3$ of the pressure, that is   $$p(t, x + a e_3) = p(t, x),$$ for all $a\in \mathbb{R}$, $t\in [0,1]$, $x\in \mathbb{R}^3$. To check that Noether's theorem yields the expected result, set $$h^\epsilon : (t, x)\in[0,1] \times \mathbb{R}^3 \to h^{\epsilon}(t,x):= x +  \epsilon e_3 \in \mathbb{R}^3.$$ By proposition 3.2. of \cite{CRUZLAS},~(\ref{invcond}) is trivially satisfied, so that, by Theorem 6.1 of \cite{CRUZLAS}, we obtain $(<v_t^{P_u}, e_z>)$, as the related $(\mathcal{F}_t^{{P_u}})-$\textit{local martingale}.  
 
\subsection{Symmetry by rotation and the kinetic momentum process} Assume the \textit{symmetry by rotation} along the axis $e_3$ of the pressure;  that is,  \begin{equation} \label{invarpress}p(t, R^\epsilon x) = p(t, x),\end{equation} for all $(t,x)\in [0,1]\times \mathbb{R}^3$ and $\epsilon \in \mathbb{R}$, where $R^\epsilon : \mathbb{R}^3\to \mathbb{R}^3$ denotes the operator of rotation, along the axis $e_3$,  with angle $\epsilon$. We consider the family of space  transformations   \begin{equation} \label{rotationdef} h^\epsilon : (t, x)\in [0,1] \times\mathbb{R}^3 \to h^{\epsilon}(t,x):= R^\epsilon  x \in \mathbb{R}^3.\end{equation} Applying Lemma 3.2. of \cite{CRUZLAS}, with $h^\epsilon$ given by~(\ref
{rotationdef}), we compute the characteristics of $\Gamma^\epsilon_\star \eta$ and we obtain  the relation $$|v^{\Gamma^\epsilon_\star \eta} \circ
\Gamma^\epsilon|_{\mathbb{R}^3}  =|v^\eta|_{\mathbb
{R}^3} \quad  \lambda \otimes \eta-a.e. ,$$ for any $\epsilon \in \mathbb{R}$. Whence, the \textit{symmetry condition}~(\ref{invcond}) is satisfied, from~(\ref{invarpress}).  Define the
stochastic process $(l_t)$, on the complete probability space $(W, \mathcal{B}(W)^{P_u},P_u)$, by  $$ l_t  :=  <W_t,e_1>_{\mathbb{R}^3} <v_t^{P_u}, e_2>_{\mathbb{R}^3}^2 - <W_t, e_2>_{\mathbb{R}^3}<v_t^{P_u}, e_1>_{\mathbb{R}^3}$$ for $t\in[0,1],$ the stochastic counterpart to the
 \textit{kinetic momentum} along $e_3$. Further denoting by $$rot \ u_t :\mathbb{R}^3 \to \mathbb{R}^3,$$ the rotational of $u(t,.)$, the weak Noether Theorem 6.1 of \cite{CRUZLAS} implies that the  corresponding process $(\mathcal{I}_t)$, defined by  $$\mathcal{I}_t:=l_t + \int_0^t <rot \ u_{1-s}(W_s),e_3> ds,$$ is the \textit{local martingale} associated to this \textit{symmetry by rotation}; otherwise stated, the stochastic kinetic momentum process $l_t$ along $e_3$  is a semi-martingale on the canonical filtered probability space of $\nu$, whose finite variation term is determined by the rotational of $u$ along the symmetry axis. The latter expresses the \textit{dissipation}, modeled through the martingale part of $P_u$, which is involved in the covariation process in the expression of $\mathcal{I}_t$.


\end{document}